\documentclass{amsart}
\usepackage[all]{xy}
\usepackage{amssymb}

\newcommand{\tr}{{\rm tr}}

\newcommand{\Div}{{\rm Div}}
\newcommand{\Id}{{\rm Id}}
\newcommand{\End}{\rm End}

\newcommand{\A}{{ \rm Aut }}

\renewcommand{\O}{{\mathcal{O}}}
\newcommand{\Z}{{\mathbb{Z}}}
\newcommand{\Q}{{\mathbb{Q}}}
\newcommand{\R}{{\mathbb{R}}}
\newcommand{\C}{{\mathbb{C}}}

\newcommand{\dv}{{\rm div}}
\newtheorem{pro}{Proposition}[section]
\newtheorem{lemma}[pro]{Lemma}

\newtheorem{theorem}[pro]{Theorem}
\newtheorem{definition}[pro]{Definition}

\begin{document}
\bibliographystyle{amsplain}

\title{Actions of Galois Groups on Invariants of number Fields}
\author{A.  Kontogeorgis}

\begin{abstract}
In this paper we investigate the connection between relations among various invariants of number fields $L^H$ 
coresponding to subgroups $H$ acting on $L$  and of linear relations among norm idempotents. 
\end{abstract}
\email{kontogar@aegean.gr}
\address{University of the \AE gean, Department of Mathematics, Karlovassi 
Samos}
\date{\today}

\maketitle

\section{Introduction}
Let $C$ be an algebraic curve defined over an algebraically closed 
field of arbitrary characteristic and 
 let $G\subset \A(C)$ be a subgroup 
of the automorphism group $G$, acting on $C$. 
  For a subgroup  $H$  of  $G$, let  
$C^H$ be the quotient group  and let $g_H$ and $\gamma_H$ be the 
the genus and the $p$-rank of the Jacobian of  
$C^H$.   In the group algebra $k[G]$ the norm idempotents $\varepsilon_H$
are defined by 
\[
\varepsilon_H=\frac{1}{|H|}\sum_{h\in H} h.  
\]
E.   Kani and M.  Rosen \cite{KaniRosen89},\cite{KaniRosen94},  studied the action of automorphisms
on the Jacobian variety of the curve, and they proved that every linear 
relation among the norm idempotents coming form subgroups $H$ of $G$ 
implies the same relations for $g_H,\gamma_H$.  
This is a generalization of results proved by R. Accola \cite{Accola70}.

They also have proved 
that this linear relation imply the same relations for the zeta 
functions of the corresponding   fields  $L^H$ where $L$ is the function field of an 
algebraic curve or a number field \cite[prop. 1.2]{KaniRosen94}, 
{\em i.e.}, 
\begin{equation}
\label{1.1}
\sum {r_H} \varepsilon_H =0 \Rightarrow \prod \zeta_H (s)^{r_H}=1.
\end{equation}
The rings of integers of 
number fields have a theory similar to that of  non singular algebraic curves, 
in the sense that the ring of integers  are Dedekind so they give rise 
to one-dimensional  affine schemes, that can be completed with the 
aid of infinite primes.   
For  number fields there is a notion of genus, 
and 
the analoga for Jacobian varieties and Tate modules can be defined.

It is known that 
a lot of information concerning a number field, can be found
in the corresponding zeta function.  
Let $L^H$ be the number field corresponding to the subgroup $H$, 
of the Galois group $G$.   
Using the characterization 
of the residues of the zeta functions for number fields at $s=1$, 
Kani and Rosen arrived at a formula, involving the class 
number $h_{L^H}$, the regulator $\mathrm{Reg}(L^H)$, and the number $w_H$ of roots 
of unity in $L^H$:
\[
\prod (h_{L^H} \mathrm{Reg}(L^H))^{r_H} =\prod w_H^{r_H},
\]
where $H$ runs over the subgroups of $G$.
The last equality was also proved by R. Brauer \cite{Brauer51} in 1950. 

In this paper we study the dependence of several group invariants
of the subfields $L^H$, corresponding to the Galois subgroups $H$, 
in terms of the linear relations among norm idempotents defined by the 
subgroup $H$.

In order to do so we give a generalization of the notion of Tate
modules for the ``Jacobian'' of a Number Field, and we consider 
the action of the Galois group on it.   
In our study, 
problems arise that are similr to those stemming from 
 wild ramification of the 
action of a group on a curve defined over a field of positive 
characteristic.   

Consider a number field $L$.   Fix a subfield $K$ such that 
the extension $L/K$ is Galois with Galois group $G$.   
For every subgroup $H$ of $G$ we define, as usual, the 
fixed field $L^H$.   The following functions from the 
set of subgroups $G$ to $\Z$ are defined:
\begin{enumerate}
\item Let $r_{H}, 2s_{H}$ be the number of real and imaginary 
embedings of $L^H$ to $\bar{\Q}$.   We set 
\[
\lambda_H=r_{H}+s_{H}-1.
\]
\item Consider the class group $Cl(L^H)$.   It is a finite 
Abelian group, hence it can be written as 
\[
Cl(L^H)=\bigoplus_{p\mid Cl(L^H)} A(H,p),
\]
where $A(H,p)$ is the $p$-part of the abelian group $Cl(L^H)$ and $A(H,p)$ in 
turn can be expressed as a finite direct sum of abelian groups $A(H,p,n)$ that are sums of 
$\lambda_{H,p,n}$ summands  of cyclic groups of order $p^n$, {\em i.e.},
\[
A(H,p)=\bigoplus _{n=1}^\infty A(H,p,n),
\]
\[
A(H,p,n)=\bigoplus_{\mu=1} ^{\lambda_{H,p,n}}  {\Z}/{p^n\Z}.
\]
Notice that in the above formulas $\lambda_{H,p,n}=0$ for all but finite integers $n$, 
and that $A(H,p,n)$ are free $\Z/p^n \Z$-modules of rank $\lambda_{H,p,n}$.
\item Consider the group $\mu(L^H)$ of units of finite order 
in the field $L^H$.   It is a cyclic group and can be written 
as 
\[
\mu(L^H)=\bigoplus_{p \mid | \mu(L^H)|} \Z/{p^{\nu(H,p)}\Z}
\]
\end{enumerate}
We will  show that 
the above functions $\lambda_H,\lambda_{H,p,n}$, behave like the 
$p$-ranks of the Jacobians of algebraic curves.   
Namely, we prove that every linear relation among norm idempotents of the 
subgroup implies the same relations for the $\lambda$ functions:
\begin{theorem}
Let $L/K$ be a Galois extension with Galois group $G$ of order $n$.  
Every relation 
\[
\sum r_H \varepsilon_H =0
\]
among norm idempotents implies relations 
\begin{align*}
\sum r_H \lambda_H &= 0 \\
\sum r_H \lambda_{H,p,n} &=0 \\
\sum r_H \nu(H,p) & =0 \mbox{ for every } p \nmid n.  
\end{align*}
\end{theorem}
On a number field $L$ we can define the notion of the Arakelov 
genus $g_L$, so it is interesting to ask whether a relation 
among norm idempotents implies the same relation among Arakelov 
genera.   The answer is yes provided we have ``tame ramification''
in the group of units that are contained in $L$, {\em i.e.},
\begin{pro}\label{pro1.2}
Let $L/K$ be a Galois extension with Galois group $G$.   
Let $w_L$ be the order of the group of units contained in $L$.   
Consider the set $S$ of subgroups $H<G$, such that $(|H|,w_L)=1$.   
Every linear relation $\sum_{H\in S} r_H \varepsilon_H=0$ among 
norm idempotents corresponding to subgroups $H\in S$, implies the 
same relation among the Arakelov genera $g_{L^H}$, of the 
fixed fields $L^H$.   In particular, if the 
order of the Galois group is prime to $w_L$, then 
$\sum r_H \varepsilon_H=0$ implies $\sum r_Hg_{L^H}=0$.   
\end{pro}
In \cite{GeerSchoof} G. Van der Geer and R.  Schoof 
introduced the notion of effectivity of an 
Arakelov divisor, a notion that is 
close to the definition  of the effectivity 
of a divisor on an algebraic curve.   This notion  gives rise to 
a new notion of $H^0(D)$, for Arakelov divisors $D$ and introduces
naturally a new  invariant  $\eta_L$ for the number field $L$:
\[
\eta_L :=
\left(\sum_{x \in \O_L} e^{-\pi || x ||^2_{L,0}} \right)
,
\]
where $||\cdot||_{L,0}$ is the metric on the Minkowski space of the number 
field $L$ defined by 
\[||x||_{L,0}^2=\sum |\sigma(x)|^2. \]
Given a relation $\sum n_H \epsilon_H=0$,
we will prove a formula for the $\eta$-invariants corresponding to 
subfields $L^H$ of $L$. 
In order to do so we have to change the model at the  infinite primes by 
considering a different metric  $||\cdot||_{L,A}$ on the Minkowski
vector space. This metric is defined in (\ref{fullmetric}).
We introduce the invariants 
\[
\eta_{A}(L) :=
\left(\sum_{x \in \O_L} e^{-\pi || x ||^2_{L,A}} \right)
,
\]
for every divisor  $A$  supported at infinite primes.  We will prove the following:
\begin{pro} \label{new-eta-rel}
Let $\sum n_H \epsilon_H=0$ be a linear relation among norm idempotents.
If $\mathbb{P}(L^H,\mathbb{R})$ (resp. $\mathbb{P}(L^H,\mathbb{C})$) 
denotes the real (resp. complex) infinite primes and 
 \[
B(H)= -\frac{\log(|H|)}{2\pi}\sum_{\sigma\in \mathbb{P}(L^H,\mathbb{R})}\sigma
 -\frac{\log(|H|/2)}{\pi}\sum_{\sigma\in \mathbb{P}(L^H,\mathbb{C})}\sigma.
\]
is a divisor supported on infinite primes of the field $L^H$, then the 
following formula holds
\[
0=\sum_H  \lambda_H \eta_{B(H)}(L^H).
\]
\end{pro}

{\bf Aknowledgement:} The author wishes to thank professor G. Van der Geer for 
his remarks and commends.
\section{Notations}
Let $K$ be a number field with ring of algebraic integers $\O_K$.   
We will follow the notation of the book of J.   Neukirch \cite{Neu}. 

An Arakelov divisor of $K$, is a formal sum 
\[
D=\sum_p \nu_p p, 
\]
where $p$ runs over the finite and infinite primes of $K$, and $\nu_p\in \Z$
if $p$ is a finite prime and $\nu_p\in \R$ if $p$ is an infinite prime.  
We will denote by 
\[
{\rm Div}(\bar{\O}_K)\cong {\rm Div}(\O_K )\times \oplus_{p\mid \infty} \R p
\]
the set of Arakelov divisors on $K$.  
There is a canonical homomorphism 
\[
\dv: K^* \rightarrow {\rm Div}(\bar{\O}_K)
\]
sending $f\in K^*$ to $\sum_p u_p(f)p$, where $v_p(f)$ is the normalized 
$p$-adic valuation of $f$ if $p$ is a finite prime, and 
$v_p(f)=-\log|\tau(f)|$, where $\tau\in {\mathrm Hom}_\Q (K,\bar{\O}_K)$ is the 
monomorphism corresponding to the infinite prime $p$.   The Arakelov 
class group $\mathrm{CH}^1(\bar{\O}_K)$ is defined as
\[
\mathrm{CH}^1(\bar{\O}_K)=\frac{{\rm Div}(\bar{\O}_K)}{\dv(K^*)}
\]
and it is equipped with the quotient topology.   Since $\prod_p |f_p|=1$
we can define on $\mathrm{CH}^1(\bar{\O}_K)$ a continuous function 
\[
\deg:\mathrm{CH}^1(\bar{\O}_K)\rightarrow \R
\]
sending $D=\sum \nu_p p $ to $\sum_p \nu_p \log(N(P))$, 
where $N(P)$ denotes the norm of $P$.   
The kernel of the degree map is a compact group denoted by 
$\mathrm{CH}^1(\bar{\O}_K)^0=:J_K$.   It can be proved \cite[Satz 1.11]{Neu} that 
$J_K$ is given by the short exact sequence:
\[
1\rightarrow H/\Gamma \rightarrow J_K \stackrel{\phi}{ \rightarrow}
Cl(K) \rightarrow 1, 
\]
where $H/\Gamma$ is homeomorphic to a torus of dimension $r+s-1$ and 
$Cl(K)$ is the ordinary class group of the number field $K$.   

Following the theory of Jacobian varieties on a curve we set 
\[
J_{K,p^n}:=\{ \mbox{Elements in } J_K \mbox{ of order } p^n.  \}
\]
where $p$ is a prime number of $\Z$.   For the $p$-part of $J_K$ we have 
the following short exact sequence:
\[
1\rightarrow \left( {\Z}/{p^n\Z} \right)^{r+s-1} \rightarrow 
J_{K,p^n} \rightarrow Cl(K)_{p^n} \rightarrow 1, 
\]
where $Cl(K)_{p^n} \cong \oplus_{i=1}^{\lambda_{K,n}} {\Z}/{p^n\Z}$
is the subgroup killed by multiplication by $p^n$. Using the 
classification theorem of finite Abelian groups we can write
\[
J_{K,p^n}\cong \bigoplus _{i=1}^{r+s-1+\lambda_{K,n}}{\Z}/{p^n\Z}. 
\]
The groups $J_{K,p^n}$ form an inverse system and we can define the 
inverse limit forming the Tate module of $J_K$ at $p$.   
Namely we set
\[
T_p(J_K)=\lim_{\leftarrow} J_{K,p^n}. 
\]
The Tate module is a free $\Z_p$-submodule of rank $s+r-1$.   Since the 
order of the ordinary class group $Cl(K)$ is finite $\lambda_{K,n}=0$, 
for large $n$, and this implies that the information of the $p$-part 
of the class group is lost after taking the inverse limit.   
We will study the $p$-part $Cl(K)_p$ of the class group 
separately.   

The action of $G$ on the primes of $L$ induces a representation
\[
\rho:G \rightarrow \End(J_L),
\]
and since  endomorphisms of $J_L$ preserve the orders of the elements
in the class group we can define representations:
\[
\rho_p:G \rightarrow \End(J_{L,p}). 
\]
Every $\rho_p$ gives rise to a representation 
\[
\hat{\rho}_p:\Q_p[G] \rightarrow \End^0(T_p(J_L)):=\End(T_p(J_L))
\otimes_{ \Z} \Q_p
\]
and to a representation 
\[
\tilde{\rho}_p: \Z_{(p)}[G] \rightarrow \End(Cl(L)_p)\otimes _{\Z} \Z_{(p)},
\]
where $\Z_{(p)}$ denotes the localization of the integer ring 
with respect to the prime ideal $p$.  
We define the $\Q_p$
vector space $V_p(J_L):=T_p(J_L)\otimes_{\Z _p} \Q_p$, so 
$\End^0(T_p(J_L))\cong \End(V_p).$
Since there is $p$-torsion on the $\Z$-module $Cl(L)_p$ we cannot 
tensor by a field, without trivializing.  
The closest structure to vector space we can 
obtain without trivializing, is by tensoring with the 
localization $\Z_{(p)}$.   
So we define the $\Z_{(p)}$-module
 $V_p(Cl(L))$ by $V_p(Cl(L)):=Cl(L)_p\otimes_\Z \Z_{(p)}$.  
The $p$-part of the class group can be factored, by the classification 
theorem of Abelian groups, as follows:
\[
Cl(L)_p =\bigoplus_{\nu=1}^\infty \bigoplus_{\mu=1}^{\lambda_{\{1\},p,\nu}}
{\Z}/{p^\nu \Z},
\]
and 
\[
V_p(Cl(L))=\bigoplus_{\nu=1}^\infty \bigoplus_{\mu=1}^{\lambda_{\{1\},p,\nu}}
{\Z_{(p)}}/{p^\nu \Z_{(p)}}.
\]
Since endomorphisms that came from $\Z_{(p)}[G]$ preserve the order 
of the group, the representation $\tilde{\rho}_p$ can be factored 
as a sum of matrix representations:
\[
\tilde{\rho}_{p,\nu}:\Z_{(p)}[G] \rightarrow M_{\lambda_{\{1\},p,\nu}}
\big(
{\Z_{(p)}}/{p^\nu \Z_{(p)}}\big), 
\]
where $M_{r}(R),$
denotes the $r\times r $ matrices with coefficients from the ring $R$.    
We define the trace of $\tilde{\rho}_p$ to be the sequence
$\tr(\tilde{\rho}_p):=(\tr(\tilde{\rho}_{p,\nu}))_\nu$.   Obviously, 
$\tr(\tilde\rho_{p,\nu})=0$ for all but finite $\nu$.  
\section{Field extensions}
Let $K,L$ be two number fields and let
\[
\tau:K \rightarrow L
\]
be a homomorphism.   For an Arakelov divisor $D=\sum_P \nu_P P$ of 
$L$ we define:
\[
\tau_*(D):=\sum_p \left( \sum_{P\mid p} \nu_P f_{P/p} \right) p \in \Div(\bar{\O}_K),
\]
where $f_{P/p}$ denotes the inertia degree of $P$ over $\tau K$ and 
$P\mid p$ means that $\tau p = P\mid_{\tau K}$.  
Conversely for an Arakelov divisor $D=\sum_p \nu_p p$ of $K$ we define
\[
\tau^*(D)=\sum_p \sum_{P\mid p} \nu_p e_{P/p}P\in \Div(\bar{\O}_L),
\]
where $e_{P/p}$ denotes the ramification index of $P$ over $\tau K$.   
The maps $\tau_*,\tau^*$ induce maps 
\[
\tau_*: \mathrm{CH}^1(\bar{\O}_L)\rightarrow \mathrm{CH}^1(\bar{\O}_K)
\]
and 
\[
\tau^*:  \mathrm{CH}^1(\bar{\O}_K)\rightarrow \mathrm{CH}^1(\bar{\O}_L)
\]
such that $\tau_* \circ \tau^* =[L:K]$ and $\deg(\tau_*D)=\deg(D)$,
$\deg(\tau^*D)=[L:K]\deg(D)$ \cite[p. 204]{Neu}.
By the above formulas for the degree, we have that there are
well defined homomorphisms
\begin{equation}
\label{3.1}
\tau_*:J_L \rightarrow J_K,
\end{equation}
and 
\begin{equation}
\label{3.2}
\tau^*:J_K \rightarrow J_L.  
\end{equation}
\begin{definition}
Let $\Z_{(p)}$ denote the localization of the ring of integers
with respect to a prime ideal.  We will denote by $R$ either a field of characteristic zero or 
$\Z_{(p)}$. 
\end{definition}
\begin{lemma} \label{le3.1}
\label{zp}
Let $V,W$ be two finitely generated $R$-modules.
  Suppose that there are two $R$-module 
homomorphisms $f_{V,W}:V \rightarrow W$, $f_{W,V}:W\rightarrow V$,
such that 
\[
f_{V,W}\circ f_{W,V}=n\; \Id_{W}
\]
and 
with $(n,p)=1$ if $R=\Z_{(p)}$.  Then there is a map:
\[
\phi:\End(W) \rightarrow \End(V),
\]
such that $\tr(a)=\tr(\phi(a))$.   In particular, if $a=\Id_W,$ then 
$\phi(a)\in \End(V),$ has trace equal to ${\rm rank}(W)$
\end{lemma}
\begin{proof}
For every $a\in \End(W)$ we define $\phi(a)\in \End(V)$ by
\[
\phi(a):=\frac{1}{n} \; f_{W,V}\circ a \circ f_{V,W}.
\]
Since $n$ is an invertible element in $R$ the map $f_{V,W}$ is onto. 
We consider the following short exact sequence of $R$-modules:
\[
\xymatrix{
0 \ar[r] & \ker f_{V,W} \ar[r] & V \ar[r]^{f_{V,W}} 
 & W \ar@/^1pc/[l]^{\frac{1}{n} f_{W,V}}  \ar[r] & 0 
}
\]
By construction, $\phi(a)$ is zero on $\ker f_{V,W}$, and 
$\tr(a)=\tr(\phi(a))$.   In particular, for $a=\Id_W$ we 
have that $\tr(\Id_W)={\rm rank}(W)$, hence $\tr(\phi(a))={\rm rank}(W)$
\end{proof}
\begin{lemma}\label{le3.2}
For a given group $G$, let $S$ be the following set of subgroups of $G$:
\[
S:=\begin{cases}
\mbox{all subgroups of } G & \mbox{if } $R$ \mbox{ is a field }\\
H<G, p\nmid |H| & \mbox{if } R \mbox{ is } \Z_{(p)}
\end{cases}
\]
Let $V$ be a free  $R$-module.   
To every $H\in S$ we attach a free  $R$-module $V(H)$ and two 
$R$-module homomorphisms $f^H: V(H)\rightarrow V$ and 
$f_H:V \rightarrow V(H)$, such that 
\[
f_H \circ f^H = |H|\cdot \Id_{V(H)}
\]
and
\[
f^H\circ f_H =\sum_{h\in H}h.   
\]
Moreover, there is a map $\phi:\mathrm{End}(V(H))\rightarrow \mathrm{End}(V)$,
such that $\phi(\mathrm{Id}_{V(H)})=\epsilon_H$ and 
\[
\mathrm{tr}(\epsilon_H)=\mathrm{rank}_R V(H).
\]
Under the above assumptions, if $\sum_{H\in S} n_H\varepsilon_H=0$
then $\sum _{H\in S} n_H {\rm rank_RV(H)}=0$.   
\end{lemma}
\begin{proof}
We apply lemma (\ref{le3.1}) for $W=V(H)$, $f_{W,V}=f^H$ and 
$f_{V,W}=f_H$.   The map 
$\phi(\Id_{V(H)})=\frac{1}{|H|}f^H\circ \Id_{V(H)} \circ f_H=\varepsilon_H$,
so $\tr(\varepsilon_H)={\rm rank}_R(V(H))$.  
\end{proof}
{\bf Remark:}
In the above lemma it was necessary to restrict ourselves 
to subgroups of order not divisible by $p$.  
The problems that appear for groups divisible by $p$, 
are of similar nature with the problems that appear 
in wild ramified extensions of rings.   
Indeed, in  the case of wild ramification of an extension 
of rings $S/R$ with Galois group $G$, the ring $S$ is not 
$R[G]$-projective.   
\begin{pro}
Let $\tau: K \hookrightarrow L$ be an inclusion 
of number fields such that $L/K$ is Galois with Galois
group $H$.   There 
is a map $\phi:\End^0(T_p(J_K))\rightarrow \End^0(T_p(J_L))$
such that 
\begin{itemize}
\item  
$\hat{\rho}_p(\varepsilon_H)=\phi(\Id_{T_p(J_K)})$
\item For any $\alpha\in \End^0(T_p(J_K))$ we have
\[
\tr(\phi(\alpha)\mid V_p(J_L)=\tr(\alpha\mid V_p(J_K)))
\]
\end{itemize}
\end{pro}
\begin{proof}
The homomorphisms defined in  (\ref{3.1}),(\ref{3.2}) can be extended
linearly to maps $V(\tau_*)$,$V(\tau^*)$ from 
$V_p(J_K)$ to $V_p(J_L)$. 
The first assertion is a direct application of lemma \ref{le3.2} if we set
$V(H)=V(T_p(J_K)))$, $V=V(T_p(J_L))$, $f_H=V(\tau_*)$, $f^H=V(\tau^*)$.
  The second assertion is an application of \ref{le3.1}.  
\end{proof}
\begin{pro}
Let $\tau: K \hookrightarrow L$ be an inclusion 
of number fields such that $L/K$ is Galois with Galois
group $H$, of order $n$.   For every prime $p$, $p\nmid n$
there is a map $\phi_2:\End^0(Cl_{K,p})\rightarrow \End^0(Cl_{L,p})$
such that 
\begin{itemize}
\item
$\tilde{\rho}_p(\varepsilon_H)=\phi_2(\Id_{Cl_{K,p}})$
\item For every $\alpha \in \End^0(Cl_{K,p})$ we have 
\[
\tr(\tau^*(\alpha)\mid V_p(Cl_L))=\tr(\alpha\mid V_p(Cl_K))
\]
\end{itemize}
\end{pro}
\begin{proof}
As before, the homomorphisms defined in (\ref{3.1}),(\ref{3.2})
can be extended $\Z_{(p)}$-linearly to maps 
$V(\tau^*)$,$V(\tau_*)$ from $V_{p,\nu}(Cl_K)$ to $V_{p,\nu}(Cl_L)$, 
where
\[
V_p(Cl_K)=\bigoplus_{\nu=1}^\infty V_{p,\nu}(Cl_K),
\]
and 
$V_{p,\nu}(Cl_K):=\oplus_{\mu=1}^{\lambda_{H,p,\nu}}
\frac{\Z_{(p)}}{p^\nu \Z_{(p)}}$.   The desired result is a direct 
consequence of lemma \ref{le3.2} if we set 
$V(H)=V_p(Cl_K)$, $V=V_p(Cl_L)$, $f_H=V(\tau_*)$, $f^H=V(\tau^*)$. 
\end{proof}
\begin{pro}
For every Galois extension $L/L^H$, with Galois group $H$, we have
\[
\tr(\hat{\rho}_p(\varepsilon_H))=\lambda_H.   
\]
For all $p$ that do not divide the order of $H$ we have
\[
\tr(\tilde{\rho}_p(\varepsilon_H))=\lambda_{H,p,\nu}.
\]
\end{pro}
\begin{proof}
This result is clear for $H=\{\Id\}$.   By the above two propositions
it is also clear for the general $H$,since
\[
\mathrm{tr}(\hat{\rho}_p(\varepsilon_H))=\mathrm{tr}(\phi\big(\mathrm{Id}_{T_p(J_K)}\big)=\lambda_H
\]
and
\[
\tr(\tilde{\rho}_p(\varepsilon_H))=
\mathrm{tr}\big(\phi_2(\mathrm{Id}_{V_p(Cl_{K})}\big)=\lambda_{H,p,\nu}.
\]
\end{proof}
We will now study the action of the Galois group on the 
group $\mu(L),\mu(L^H)$ of units contained in the fields 
$L,L^H$ respectively.   
\begin{pro}\label{3.7}
Let $L/L^G$ be a Galois extension with Galois group $G$.  Let 
$w_L$ be the order of the group $\mu(L)$ of units of finite order.   
For a fixed prime $p\mid w_L$ we define the set $S_p$ of subgroups 
$H$ of $G$ with the property 
$H\in S_p$ if and only if $p\nmid |H|$.   Let $\nu(H,p)$ be the 
valuation at $p$ of the order $w_{L^H}$ of the unit group of $L^H$.   
Every linear relation of the form $\sum_{H\in S_p} n_H \varepsilon_H=0$, 
among norm idempotents of groups $H\in S_p$, implies the same relation 
$\sum n_H \nu(H,p)=0$.   
\end{pro}
\begin{proof}
Consider the norm 
\[
N_{L/L^H}:\mu(L) \rightarrow \mu(L^H),
\]
and the inclusion function $i_{L^H,L}:\mu(L^H)\rightarrow \mu(L)$.   
We have
\[
N_{L/L^H}\circ i_{L^H,L}=|H| \cdot \Id_{\mu(L^H)}
\]
and 
\[
i_{L^H,L}\circ N_{L/L^H}=\sum_{h\in H}h.  
\]
The unit group $\mu(L)$ is a cyclic group of order $m$, and 
$\mu(L^H)$ is a subgroup of the cyclic group $\mu(L)$.   The group 
$\mu(L)$ can be considered as a direct sum 
\[
\mu(L)=\bigoplus_{i=1}^r \frac{\Z}{p_i^{\nu(\{1\},p_i)}}=\bigoplus_{i=1}^r
\mu_i(L),
\]
where $p_i$ are the different prime divisors of $m$.   Each direct summand
$\mu_i(L)$, gives rise to a $\Z_{(p_i)}$-module $\mu_i(L)\otimes_\Z \Z_{(p_i)}$.
We can have a similar decomposition of $\mu(L^H)$ as a direct sum of $ \Z_{(p_i)}$-modules:
\[
\mu(L^H)=\bigoplus_{i=1}^r \frac{\Z}{p_i^{\nu(H,p_i)}}=\bigoplus_{i=1}^r
\mu_i(L^H),
\]  
The $N_{L/L^H}$ and 
$i_{L,L^H}$ group homomorphism give rise to $\Z_{(p_i)}$-module 
homomorphisms from $\mu_i(L)\otimes _\Z \Z_{(p_i)}$ to $\mu_i(L^H)\otimes _\Z \Z_{(p_i)}$.   
The desired result follows by lemma \ref{le3.2}.
\end{proof}
\section{Analytic methods}
In this section we give an analytic proof of proposition \ref{pro1.2}.
Consider the zeta function of an algebraic number field $L$,
\[
\zeta_L(s):=\sum_{A\in I_L}\frac{1}{N(A)^s},\;\;Re(s)>1,
\]
where $A$ runs over the integral ideals $I_L$ of the ring of 
integers of $L$. It is known that $\zeta_L(s)$ admits a meromorphic 
extension in $\C\setminus\{1\}$ with only one  pole at $s=1$.  
Moreover the residue at $s=1$ can be computed \cite[Satz VII 5.11]{Neu}
\[
\mathrm{Res}_{s=1}\zeta_L(s)=\lim_{s\rightarrow 1^+} (s-1)\zeta_L(s)=
\frac{2^r(2\pi)^s \mathrm{Reg}(L)}{|\mu(L)|\sqrt{|D_L|}}h_L.  
\]
The Arakelov genus $g_L$ of the number field $L$
is defined by 
\[
g_L=\log\frac{|\mu(L)|\sqrt{|D_L|}}{2^r(2\pi)^s},
\]
therefore 
\[
\mathrm{Res}_{s=1}\zeta_L(s)=e^{-g_L} \mathrm{Reg}(L) h_L.  
\]
Let $\sum n_H \varepsilon_H $ be a norm idempotent relation.   
The product formula (\ref{1.1}) implies that
\[
\lim_{s\rightarrow 1^+} \sum n_H \log((s-1)\zeta_{L^H}(s))=\sum n_H 
\lim_{s \rightarrow 1^+} \log(s-1).
\]
The left hand side is finite (the regulator of every number field 
is not zero), therefore $\sum n_H=0$ and moreover
\[
\sum n_H \left(  -g_{L^H}+ \log(\mathrm{Reg}(L^H) h_{L^H}) \right) =0 \Rightarrow
\]
\begin{equation}
\label{4.1}
\sum n_H g_{L^H} = \sum n_h \log(\mathrm{Reg}(L^H) h_{L^H}).
\end{equation}
{\bf Remark:} The relation $\sum n_H =0$ can also be proved by applying the 
character of the trivial representation on the sum $\sum n_H \varepsilon_H$.  

On the other hand using the analytic continuation of the $\zeta_L(s)$ we 
can prove that 
\[
\lim_{s\rightarrow 0} \zeta_L(s)/s^{r+s-1}=-\frac{h_L \mathrm{Reg}(L^H)}{|\mu(L)|},
\]
therefore
\begin{equation}
\label{4.2}
\sum n_H \log(h_{L^H} \mathrm{Reg}(L^H))=\sum n_H \log |\mu(L^H)|. 
\end{equation}
Combining (\ref{4.1}),(\ref{4.2}), we arrive at 
\begin{equation}
\label{4.3}
\sum n_H g_{L^H}=\sum n_H \log(\mu(L^H)).
\end{equation}
For every $H$, such that $(|H|,|\mu(L)|)=1$, we write 
$|\mu(H)|=\prod_{i=1}^r p_i^{\nu_{H,p_i}}$.   Therefore the right hand 
side of (\ref{4.3}) is written
\[
\sum_H n_H \log(|\mu(L^H)|)=\sum_H n_H \log(\prod_{i=1}^r p_i^{\nu_{H,p_i}})=
\sum_{i=1}^r \log(p_i)\sum_H n_H \nu_{H,p_i}=0
\]
by proposition \ref{3.7} and the proof of proposition \ref{pro1.2} is now complete. 

\section{The $\eta$ invariant}

In order to apply the proof given in previous sections  we would like 
to realize the function 
\[
\sum_{x \in \O_L} e^{-\pi || x ||^2_{L}}
\]
as the trace of a suitable linear operator. 

If $L$ is a number field we will denote by $r_L$ the number of real embedings and 
by $s_L$ the number of complex nonequivalent embeddings. The Minkowski space 
is defined by 
\[
M(L)=\mathbb{R}^{r_L} \times \mathbb{C}^{s_L}.
\]
We will define the set of real infinite primes of $L$ by $\mathbb{P}(L,\mathbb{R})$ 
and by $\mathbb{P}(L,\mathbb{C})$ the set of complex infinite primes.
The field $L$ can be embedded on the Minkowski space by the map
\[
i_L: L \rightarrow M(L),
\]
\[
x \mapsto (\sigma_1(x),\ldots, \sigma_{r_L}, \sigma_{r_L+1},\ldots, \sigma_{r_L+s_L}).
\]
Every divisor 
\[
D=\sum_{\sigma \in \mathbb{P}(L,\mathbb{R})} a_\sigma \sigma + 
\sum_{\sigma \in \mathbb{P}(L,\mathbb{C})} a_\sigma \sigma,
\]
 supported on the set of infinite primes gives rise to the metric 
\begin{equation} \label{fullmetric}
||x||_{L,D}^2=
\sum_{\sigma \in \mathbb{P}(L,\mathbb{R})}  |x_\sigma|^2 e^{- 2 a_\sigma}  + 
\sum_{\sigma \in \mathbb{P}(L,\mathbb{C})}  |x_\sigma|^2  2 e^{-a_\sigma} ,
\end{equation}
where $x=(x_\sigma)$ is an element of the Minkowski space $M(L)$.

Let $L/K$ be a Galois extension of number fields with Galois group $Gal(L/K)=H$. 
An infinite  complex  prime $\sigma$ of $K$ is extended to $|H|$ infinite primes of $L$. 
Moreover, 
\begin{equation} \label{complex-ramif}
\sigma=\sum_{i=1}^{|H|} \sigma_i
\end{equation}
On the other hand an infinite real prime $\tau$ of  $K$ gives rise to $a(\tau)$ real infinite primes 
$\{\sigma_1,\ldots,\sigma_a\}$ 
of $L$ and $b(\tau)$ pairs 
$\{
\sigma_{a(\tau)+1},\ldots,\sigma_{a(\tau)+b(\tau)},\overline{\sigma_{a(\tau)+1}},\ldots,\overline{\sigma_{a(\tau)+b(\tau)}}
\}$
  complex infinite primes of $L$, where  $a(\tau)+2b(\tau)=|H|$.
So the real infinite prime $\sigma$ is decomposed in $L$ as follows: 
\begin{equation} \label{real-ramif}
\sigma=\sum_{i=1}^{a(\tau)} \sigma_i + \sum_{j=1}^{b(\tau)} \sigma_{a(\tau)+j}^2.
\end{equation}
\begin{lemma}
Let $L/K$ be a Galois extension of number fields, with Galois group $H$.
Consider the set $\mathbb{P}(L,\infty)$ (resp. $\mathbb{P}(K,\infty)$) of infinite 
primes of $L$ (resp. $K$) and let $r_L$, $s_L$ (resp. $r_K$,$s_K$) denote the 
number of real and complex embedings of $K$ (resp. $L$). 
Let $D$ be a divisor supported at the infinite primes of $L$, such that $D$ is $H$-invariant
{\em i.e.},
\[
D=\sum_{\sigma \in \mathbb{P}(L,\infty)} a_\sigma \sigma = \sum_{\tau \in \mathbb{P}(K,\infty)} a_\tau \sum_{\sigma \mid \tau}  \sigma.
\]
Let us denote by $D^H$ the divisor 
\[
D^H=\sum_{\tau \in \mathbb{P}(K,\infty)} a_\tau \tau.
\]
If $||\cdot||_{L,D}$ is the metric on the Minkowski space 
$\mathbb{R}^{r_L}\times \mathbb{C}^{s_L}$ introduced by $D$ and $||\cdot||_{K,D^H}$ is the metric on the Minkowski 
space  $\mathbb{R}^{r_K}\times \mathbb{C}^{s_K}$ introduced by $D^H$,  
then for every $x\in K\subset L$  considered as an element on the spaces $\mathbb{R}^{r_L}\times \mathbb{C}^{s_L}$ and 
$\mathbb{R}^{r_K}\times \mathbb{C}^{s_K}$ we have:
\begin{equation} \label{change-metric}
||i_L(x)||^2_{L,D}=|H|\cdot || i_K(x)||^2_{K,D^H}.
\end{equation} 
\end{lemma}
\begin{proof}
Let $x\in K \subset L$. 
We compute:
\[
||i_K(x)||_{K,D^H}^2=\sum_{\tau \in \mathbb{P}(K,\mathbb{R})}  |\tau(x)|^2 e ^{-2 a_\tau} +
\sum_{\tau \in \mathbb{P}(K,\mathbb{C})}  |\tau(x)|^2 2 e ^{- a_\tau}.
\]
If $\tau \in \mathbb{P}(K,\mathbb{R})$ then the contribution to the norm of the infinite primes 
 $ a_\tau\sum_{\sigma \mid \tau} \sigma$ above $\tau$ is according to (\ref{real-ramif}):
\[
|\tau(x)|^2
 \left(
a(\tau) e^{-2 a_\tau} + b(\tau) 2 e^{-2 a_\tau}
\right)
=|\tau(x)|^2 |H|e^{-2 a_\tau}.
\] 
On the other hand if $\tau$ is a complex prime, then the contribution to the norm of the 
infinite primes above $\tau$ is according to (\ref{complex-ramif})
\[
|\tau(x)|^2 |H| 2 e^{- a_\tau}.
\]
The desired result follows by adding all the contributions of primes of $L$ above each 
infinite prime $\tau$ of $K$.
\end{proof}
To every number field $L$ we attach the 
 Hilbert space $V_L$  consisting of functions $f:\O_L \rightarrow \mathbb{R}$, 
such that $\sum_{y \in \O_L} |f(y)|^2 < \infty$.

Let $H$ be a group acting on the number field $L$. The space $V_L$ is acted on 
by $H$ as follows:
\[
f^h(x)=f(hx), \mbox{ for } h \in H.
\]
Let $V_{L^H}$ be the Hilbert space of functions $f:\O_{L^H} \rightarrow \mathbb{R}$ such that 
$\sum_{y\in \O_{L^H}} | f(y)|^2 < \infty$.

The norm idempotent $\epsilon_H$ induces a map $\epsilon_H^*:V_{L^H} \rightarrow V_{L}$, 
sending the function $f:\O_{L^H}\rightarrow \mathbb{R}$ to the function 
$f\circ \epsilon_H: \O_L \rightarrow \mathbb{R}$. Moreover we will consider the 
restriction map $rest: V_{\O_L} \rightarrow V_{\O_{L^H}}$ sending a 
function $f:\O_L \rightarrow \mathbb{R}$, to the restriction on $\O_{L^H}$. 

Since the vector spaces we treat are of infinite dimension we can not use the 
trace of the identity map. Instead, we consider the diagonal linear operator 
\[
T_D: V_{\O_{L}} \rightarrow V_{\O_L},
\]
sending a function 
\[
\O_L \ni x \mapsto f(x) \mbox{ to } \O_L \ni x \mapsto T_Df(x)=e^{-\pi ||x||^2_{L,D}} f(x).
\]
Let $\delta_x(\cdot)$ denote the basis functions
\[
\delta_x(y)=\left\{ 
\begin{array}{l}
1 \mbox{ if } x=y \\
0 \mbox{ if } x\neq y
\end{array}
\right.
\]
We observe that the trace of the linear operator $T\circ \epsilon $ is 
given by 
\begin{equation} \label{trTH}
{\mathrm tr}(T\circ \epsilon_H)
=\sum_{x\in \O_L} \langle T\circ \epsilon_H^*(\delta_x),\delta_x \rangle
=\sum_{x\in \O_{L^H}} e^{-\pi ||x||^2_{L,D}}.
\end{equation}
Indeed, if $x\in \O_{L}$ is an element of $\O_{L^H}$ then 
$\epsilon_H(x)=x$, and if $x\in \O_L \backslash \O_{L^H}$ 
then $\epsilon_H(x) \neq x$ since $\epsilon_H(x) \in \O_{L^H}$.
We compute 
\[
\langle \epsilon_H^*(\delta_x),\delta_x \rangle= \left\{
\begin{array}{l}
1 \mbox{ if } x \in \O_{L^H} \\
0 \mbox{ if } x \in \O_L \backslash \O_{L^H}
\end{array}
\right.,
\]
and the formula (\ref{trTH}) follows.
Suppose now that $D$ is an $H$-invariant divisor. 
Then equation (\ref{change-metric}) together with (\ref{trTH}) gives that 
\begin{equation} \label{second-trace}
{\mathrm tr}(T\circ\epsilon_H)=\sum_{x\in \O_{L^H}} e^{-\pi || x ||^2_{K,D^H} |H|}.
\end{equation}
\begin{pro}
Given a number field $K$  and a divisor $A$ supported on infinite primes, define the numbers
\[
\eta_A(K)=  \sum_{x\in \O_K} e^{-\pi || x ||^2_{K,A}}.
\]
If $\sum_H \lambda_H \epsilon_H$ is a linear relation among norm idempotents,
and $D$ is an $H$-invariant divisor supported on infinite primes of $L$ 
then the following relation holds
\[
0=\sum_H  \lambda_H {\eta_{D+B(H)}(L^H)},
\]
where
\[
B(H)= -\frac{\log(|H|)}{2\pi}\sum_{\sigma\in \mathbb{P}(L^H,\mathbb{R})}\sigma
 -\frac{\log(|H|/2)}{\pi}\sum_{\sigma\in \mathbb{P}(L^H,\mathbb{C})}\sigma.
\]
In particular if $D=0$ then 
\[
0=\sum_H  \lambda_H {\eta_{B(H)}(L^H)}
\]
\end{pro}
\begin{proof}
The desired result follows by linearity of the trace map composed by $T$, 
the relation $\sum_{H} \lambda_H \epsilon_H$ and 
equation (\ref{second-trace}).
\end{proof}


\begin{thebibliography}{1}

\bibitem{Accola70}
Robert D.~M. Accola, \emph{Two theorems on {R}iemann surfaces with noncyclic
  automorphism groups.}, Proc. Amer. Math. Soc. \textbf{25} (1970), 598--602.
  \MR{MR0259105 (41 \#3747)}

\bibitem{Brauer51}
Richard Brauer, \emph{Beziehungen zwischen {K}lassenzahlen von {T}eilk\"orpern
  eines galoisschen {K}\"orpers}, Math. Nachr. \textbf{4} (1951), 158--174.
  \MR{MR0039760 (12,593b)}

\bibitem{KaniRosen89}
E.~Kani and M.~Rosen, \emph{Idempotent relations and factors of {J}acobians},
  Math. Ann. \textbf{284} (1989), no.~2, 307--327. \MR{MR1000113 (90h:14057)}

\bibitem{KaniRosen94}
Ernst Kani and Michael Rosen, \emph{Idempotent relations among arithmetic
  invariants attached to number fields and algebraic varieties}, J. Number
  Theory \textbf{46} (1994), no.~2, 230--254. \MR{MR1269254 (95c:11080)}

\bibitem{Neu}
J\"urgen Neukirch, \emph{Algebraische {Z}ahlentheorie.}, Berlin etc.:
  Springer-Verlag. xiii, 595 S., 1992 (German).

\bibitem{GeerSchoof}
Gerard van~der Geer and Ren{\'e} Schoof, \emph{Effectivity of {A}rakelov
  divisors and the theta divisor of a number field}, Selecta Math. (N.S.)
  \textbf{6} (2000), no.~4, 377--398. \MR{MR1847381 (2002e:11157)}

\end{thebibliography}

\providecommand{\bysame}{\leavevmode\hbox to3em{\hrulefill}\thinspace}
\providecommand{\MR}{\relax\ifhmode\unskip\space\fi MR }
\providecommand{\MRhref}[2]{%
  \href{http://www.ams.org/mathscinet-getitem?mr=#1}{#2}
}
\providecommand{\href}[2]{#2}

\end{document}